\documentclass[a4paper,11pt,onecolumn,twoside]{article}
\usepackage{mathrsfs}
\usepackage{mathrsfs}
\usepackage{amsfonts}
\usepackage{fancyhdr}
\usepackage{amsmath,amsfonts,amssymb,graphicx,amsthm}
\begin{document}

\title{\bf Time periodic solutions of compressible fluid models of Korteweg type}
\author{{\bf Zhengzheng Chen}\\
School of Mathematics and Statistics\\
Wuhan University, Wuhan 430072, China\\[2mm]
{\bf Qinghua Xiao}\\
School of Mathematics and Statistics\\
Wuhan University, Wuhan 430072, China\\[2mm]
{\bf Huijiang Zhao}\thanks{Corresponding author.E-mail:
hhjjzhao@hotmail.com}\\
School of Mathematics and Statistics\\
Wuhan University, Wuhan 430072, China}
\date{}

\vskip 0.2cm

\maketitle

\vskip 0.2cm \arraycolsep1.5pt
\newtheorem{Lemma}{Lemma}[section]
\newtheorem{Theorem}{Theorem}[section]
\newtheorem{Definition}{Definition}[section]
\newtheorem{Proposition}{Proposition}[section]
\newtheorem{Remark}{Remark}[section]
\newtheorem{Corollary}{Corollary}[section]

\begin{abstract}
This paper is concerned with the existence, uniqueness and time-asymptotic stability of time periodic solutions to the compressible Navier-Stokes-Korteweg system effected by a time periodic external force in $\mathbb{R}^n$. Our analysis is based on a
combination of the energy method and the time decay estimates of solutions to the linearized system.

 \bigbreak
\noindent

{\bf \normalsize Keywords }  {Navier-Stokes-Korteweg
system;\,\,Capillary fluids;\,\, Time periodic solution; \,\,Energy
estimates;}\bigbreak
 \noindent{\bf AMS Subject Classifications 2010:} 35M10, 35Q35, 35B10.

\end{abstract}

\section{Introduction }
\setcounter{equation}{0}

The compressible Navier-Stokes-Korteweg system for the density $\rho>0$ and
  velocity $u=(u_1,u_2, \cdots,u_n)\\\in\mathbb{R}^n$ is written as :
  \begin{eqnarray}\label{1.1}
\left\{\begin{array}{ll}
          \rho_t+\nabla\cdot(\rho u)=0,\\[2mm]
          (\rho u)_t+\nabla\cdot(\rho u\bigotimes u)+\nabla P(\rho)-\mu\Delta u-(\nu+\mu)\nabla(\nabla\cdot u)=\kappa\rho\nabla\Delta\rho+\rho
          f(t, x).
 \end{array}\right.
\end{eqnarray}
Here, $(t, x)\in\mathbb{R}^+\times\mathbb{R}^n$, $P=P(\rho)$ is the pressure,
$\mu, \nu$ are the viscosity coefficients, $\kappa$ is the
capillary coefficient, and $f(t,x)=(f_1, f_2, f_3)(t,x)$ is a given external force. System
(1.1) can be used to describe the motion of the compressible isothermal fluids with capillarity effect of materials, see
the pioneering work by Dunn and Serrin \cite{J. E. Dunn-J.
Serrin-1985}, and also \cite{D. M. Anderson-G. B. McFadden-G. B.
Wheeler-1998, J. W. Cahn-J. E. Hilliard-1998, M. E. Gurtin-D.
Polignone-J. Vinals-1996}.

In this paper, we consider the problem (1.1) for $(\rho,u)$ around a
constant state $(\rho_\infty,0)$ for $n\geq5$, where $\rho_\infty$ is
a positive constant. Throughout this paper, we make the following
basic assumptions:

 (H1): $\mu$, $\nu$ and $\kappa$ are positive constants and
 satisfying $\nu+\frac{2}{n}\mu\geq0$.

 (H2): $P(\rho)$ is smooth in a neighborhood of $\rho_\infty$
 satisfying $P^\prime(\rho_\infty)>0$.

 (H3): $f$ is time periodic with period $T>0$.

The main purpose of this paper is to show that the problem (1.1)
admits a time periodic solution around the constant state
$(\rho_\infty, 0)$ which has the same period as $f$. By combining the
energy method and the optimal decay estimates of solutions to the
linearized system, we prove the existence of a time periodic
solution in some suitable function space. Notice that some similar
results have been obtained for the compressible Navier-Stokes
equations and Boltzmann equation, cf. \cite{H. F. Ma-S. Ukai-T.
Yang-2010, S. Ukai-2006, S. Ukai-T. Yang-2006, R. J. Duan-S. Ukai-T.
Yang-H. J. Zhao-2008}.

Precisely, Let $N\geq n+2$ be a positive integer, define the
solution space by
\begin{equation}\label{1.2}
X_M(0, T)=\left\{(\rho, u)(t,x)\left|
\begin{array}{c}
\rho(t,x)\in C(0, T; H^N(\mathbb{R}^n))\cap C^1(0, T; H^{N-2}(\mathbb{R}^n)),\\[2mm]
u(t,x)\in C(0, T; H^{N-1}(\mathbb{R}^n))\cap C^1(0, T; H^{N-3}(\mathbb{R}^n)),\\[2mm]
\nabla\rho(t,x)\in L^2(0,T; H^{N+1}(\mathbb{R}^n)),\\[2mm] \nabla
u(t,x)\in L^2(0,T; H^{N}(\mathbb{R}^n)), |||(\rho,u)|||\leq M,
\end{array}
\right.\right\}
\end{equation}
for some positive constant $M$ and with the norm
\begin{equation}\label{1.3}
|||(\rho,u)|||^2=\sup_{0\leq t\leq
T}\left\{\|\rho(t)\|_N^2+\|u(t)\|_{N-1}^2\right\}+\int_0^T\left(\|\nabla\rho(t)\|_{N+1}^2+\|\nabla
u(t)\|_{N}^2\right)dt.
\end{equation}
Then the existence of the time periodic solution can be
stated as follows.
\begin{Theorem}
Let $n\geq5, N\geq n+2$. Assume the assumptions (H1)-(H3) hold, and
$f(t, x)\in C(0,T; H^{N-1}(\mathbb{R}^n)\cap L^1(\mathbb{R}^n))$.
Then there exists a small constant $\delta_0>0$ and a constant $M_0>0$ which are dependent on $\rho_\infty$,  such that if
\begin{equation}\label{1.4}
\sup_{0\leq t\leq T}\|f(t)\|_{H^{N-1}\cap L^1}\leq\delta_0,
\end{equation}
then the problem (1.1) admits a time periodic solution $(\rho^{per}, u^{per})$ with period $T$,
satisfying
\[
(\rho^{per}-\rho_\infty,u^{per})\in X_{M_0}(0,T)
\]
Furthermore the periodic solution is unique in the following sense: if there is another time periodic solution $(\rho_1^{per}, u_1^{per})$ satisfying (1.1) with the same $f$, and
$(\rho_1^{per}-\rho_\infty,u_1^{per})\in X_{M_0}(0,T)$,
then $(\rho_1^{per}, u_1^{per})=(\rho^{per}, u^{per})$.
\end{Theorem}

To study the stability of the time periodic solution $(\rho^{per},
u^{per})$ obtained in Theorem 1.1, we consider the problem (1.1)
with the following initial date
\begin{equation}\label{1.5}
(\rho, u)(t,x)|_{t=0}=(\rho_0, u_0)(x)\rightarrow(\rho_\infty, 0),\quad as\,\, |x|\rightarrow\infty.
\end{equation}
Here $\rho_0(x)$ and $u_0(x)$ is a small perturbation of the time
periodic solution $(\rho^{per}, u^{per})$. And we have the following stability result.
\begin{Theorem}
Under the assumptions of Theorem 1.1, let $(\rho^{per}, u^{per})$ be
the time periodic solution thus obtained. If the initial date
$(\rho_0, u_0)$ be such that $\|(\rho_0-\rho^{per}(0),
u_0-u^{per}(0)\|_{N-1}$ is sufficiently small, then the Cauchy
problem (1.1), (\ref{1.5}) has a unique classical solution $(\rho, u)$
globally in time, which satisfies
\begin{equation}\label{1.6}
\begin{array}{rl}
&\rho-\rho^{per}\in C(0, \infty; H^{N-1}(\mathbb{R}^n))\cap C^1(0,
\infty; H^{N-3}(\mathbb{R}^n)),\\[2mm]
 &u-u^{per}\in C(0, \infty;
H^{N-2}(\mathbb{R}^n))\cap C^1(0, \infty; H^{N-4}(\mathbb{R}^n)).
\end{array}
\end{equation}
Moreover, there exists a constant $C_0>0$ such that
\begin{equation}\label{1.7}
\begin{array}{rl}
 &\|(\rho-\rho^{per})(t)\|^2_{N-1}+\|(u-u^{per})(t)\|^2_{N-2}+\displaystyle\int_0^t\left(\|\nabla(\rho-\rho^{per})(\tau)\|^2_{N-1}+\|\nabla(u-u^{per})(\tau)\|^2_{N-2}\right)d\tau\\[3mm]
&\leq
C_0\left(\|\rho_0-\rho^{per}(0)\|^2_{N-1}+\|u_0-u^{per}(0)\|^2_{N-2}\right),
\end{array}
\end{equation}
for any $t\geq0$ and
\begin{eqnarray}\label{1.8}
\|(\rho-\rho^{per}, u-u^{per})\|_{L^\infty}\rightarrow0 \,\,\,as\,\,
t\rightarrow\infty.
\end{eqnarray}
\end{Theorem}

Now we outline the main ingredients used in proving of our main results. For
the proof of Theorem 1.1, thanks to the time decay estimates of
solutions to the linear system (\ref{2.6}) (see Lemma 2.1 below), we can show the integral in (\ref{4.5}) is convergent.
Based on this and the elaborate energy estimates given
in Section 3, we prove the existence of time periodic solution by
the contraction mapping principle. Here, similar to the case of compressible
Navier-Stokes equations, Theorem 1.1 is obtained only in the case $n\geq5$
because of the convergence of the integral in (\ref{4.5}). Thus, how
to deal with the case $n<5$, especially, the physical case $n=3$, is
still an open problem. Theorem 1.2 is established by the energy
method. The key ingredient in the proof of Theorem 1.2,  among other
things, is to get the a priori estimates, which can be done similarly
to the estimates in Section 3.

There have been a lot of studies on the mathematical theory of the
compressible Navier-Stokes-Korteweg system. For example, Hattori and
Li \cite{H. Hattori-D. Li-1994, H. Hattori-D. Li-1996} proved the local existence and
the global existence of smooth solutions in Sobolev space. Danchin and Desjardins \cite{R. Danchin-B. Desjardins-2001} studied
the existence of suitably smooth solutions in critical Besov space. Bresch, Desjardins and Lin
\cite{D. Bresch-B. Desjardins-C. K. Lin-2003} considered the
global existence of weak solution,  then Haspot improved their
results in \cite{B. Haspot-2009}. The local existence of strong
solutions was proven in \cite{M. Kotschote-2008}. Recently, Wang and
Tan \cite{Y. J. Wang-Z. Tan-2011} established the optimal decay
rates of global smooth solutions without external force. Li \cite{Y.
P. Li} discussed the global existence and optimal $L^2$-decay rate
of smooth solutions with potential external force.

The rest of the paper is organized as follows. In Section 2, we will
reformulate the problem and give some preliminaries for later use.
In Section 3, we give the energy estimates on the linearized system
(\ref{2.3}). The proof of Theorem 1.1 is given in Section 4. In
the last section, we will study the stability of the time periodic
solution.

 {\bf Notations:} Throughout this paper, for simplicity, we
will omit the variables $t, x$ of functions if it does not cauchy
any confusion. $C$ denotes a generic positive constant which may vary in
different estimates. $\langle \cdot,\cdot\rangle$ is the inner product in $L^2(\mathbb{R}^n)$. The norm in the
usual Sobolev Space $H^s(\mathbb{R}^n)$ are denoted by $\|\cdot\|_s$
for $s\geq0$. When s=0, we will simply use
$\|\cdot\|$. Moreover, we denote $\|\cdot\|_{H^s}+\|\cdot\|_{L^1}$ by $\|\cdot\|_{H^s\cap L^1}$. If $g=(g_1, g_2, \cdots, g_n)$, then
$\|g\|=\displaystyle\sum_{k=1}^{n}(\|g_k\|^2)^{\frac{1}{2}}$.
$\nabla=(\partial_1, \partial_2, \cdots, \partial_n)$ with
$\partial_i=\partial_{x_i}, i=1,2, \cdots, n$ and for any integer $l\geq0$, $\nabla^lg$ denotes all $x$ derivatives
of order $l$ of the function $g$.
Finally, for multi-index $\alpha=(\alpha_1, \alpha_2, \cdots,\alpha_n)$, it is standard that
\[
\partial_x^\alpha=\partial_{x_1}^{\alpha_1}\partial_{x_2}^{\alpha_2}\cdots\partial_{x_n}^{\alpha_n},
\quad |\alpha|=\sum_{i=1}^{n}\alpha_i.\]

\section{Reformulated system and preliminaries}
\setcounter{equation}{0} We reformulate the system (1.1) in this
section. Firstly, set
\[ \gamma=\sqrt{P^\prime(\rho_\infty)},\quad
\kappa^\prime=\frac{\rho_\infty}{\gamma}\kappa,\quad
\mu^\prime=\frac{\mu}{\rho_\infty},\quad
\nu^\prime=\frac{\nu+\mu}{\rho_\infty},\quad
\lambda_1=\frac{\gamma}{\rho_\infty},\quad
\lambda_2=\frac{\rho_\infty}{\gamma},\] and define the new variables
\[\sigma=\rho-\rho_\infty,\quad v=\lambda_2 u,\]
then the system (1.1) is reformulated as
 \begin{eqnarray}\label{2.1}
\left\{\begin{array}{ll}
          \sigma_t+\gamma\nabla\cdot v=G_1(\sigma, v), \\[2mm]
          v_t-\mu^\prime\Delta v-\nu^\prime\nabla(\nabla\cdot
          v)+\gamma\nabla\sigma-\kappa^\prime\nabla\Delta\sigma=G_2(\sigma,
          v)+\lambda_2f,
 \end{array}\right.
\end{eqnarray}
where
\[ \begin{array}{rl}
&G_1(\sigma, v)=-\lambda_1\nabla\cdot(\sigma v),\\[2mm]
&G_2(\sigma,v)=\displaystyle-\frac{\sigma}{\rho_\infty(\sigma+\rho_\infty)}(\mu\Delta
  v+\nu\nabla(\nabla\cdot
  v))-\lambda_1(v\cdot\nabla)v-\lambda_2\left[\frac{P^\prime(\sigma+\rho_\infty)}{\sigma+\rho_\infty}-\frac{P^\prime(\rho_\infty)}{\rho_\infty}\right]\nabla\sigma.
\end{array}
\]
Notice that $G_1$ and $G_2$ have the following properties:
\begin{equation}\label{2.2}
\begin{array}{rl} &G_1(\sigma, v)\thicksim\nabla\sigma\cdot
v+\sigma \nabla\cdot
v,\\[2mm]
&G_2(\sigma,v)\thicksim \sigma\Delta v+\sigma\nabla(\nabla\cdot
v)+(v\cdot\nabla)v+\sigma\nabla\sigma.
\end{array}
\end{equation}
Here $\thicksim$ means that two side are of same order.

 Set
$U=(\sigma,v)$, $G=(G_1, G_2)$, $F=(0, \lambda_2f)$ and
\[
\mathbb{A}=\left(\begin{array}{ll}\quad\quad0\qquad\qquad\qquad \gamma div\\[3mm]
   \gamma\nabla-\kappa^\prime\nabla\Delta\qquad-\mu^\prime\Delta-\nu^\prime\nabla div
\end{array}\right),
\]
then the system (2.1) takes the form
\begin{equation}\label{2.2a}
U_t+\mathbb{A}U=G(U)+F.
\end{equation}
 We first consider the linearized
system of (\ref{2.1}):
\begin{eqnarray}\label{2.3}
\left\{\begin{array}{ll}
          \sigma_t+\gamma\nabla\cdot v=G_1(\tilde{U}), \\[2mm]
          v_t-\mu^\prime\Delta v-\nu^\prime\nabla(\nabla\cdot
          v)+\gamma\nabla\sigma-\kappa^\prime\nabla\Delta\sigma=G_2(\tilde{U})+\lambda_2f,
 \end{array}\right.
\end{eqnarray}
for any given functions $\tilde{U}=(\tilde{\sigma},\tilde{v})$
satisfying
\[\tilde{\sigma}\in H^{N+2}(\mathbb{R}^n),\quad\tilde{v}\in
H^{N+1}(\mathbb{R}^n).\] Notice that the system (\ref{2.3}) can be
written as
\begin{equation}\label{2.4}
U_t+\mathbb{A}U=G(\tilde{U})+F.
\end{equation}
By the Duhamel's principle, the solution to the system (\ref{2.3}) can
be written in the mild form as
\begin{equation}\label{2.5}
U(t)=\displaystyle\mathbb{S}(t,s)U(s)+\int_s^t\mathbb{S}(t,\tau)(G(\tilde{U})+F)(\tau)d\tau,\quad
t\geq s,
\end{equation}
where $\mathbb{S}(t,s)$ is the corresponding linearized solution
operator defined by
\[
\mathbb{S}(t,s)=e^{(t-s)\mathbb{A}},\quad t\geq s.
\]
Indeed, the corresponding homogeneous linear system to (\ref{2.3})
is
\begin{eqnarray}\label{2.6}
\left\{\begin{array}{ll}
          \sigma_t+\gamma\nabla\cdot v=0, \\[2mm]
          v_t-\mu^\prime\Delta v-\nu^\prime\nabla(\nabla\cdot
          v)+\gamma\nabla\sigma-\kappa^\prime\nabla\Delta\sigma=0,\\[2mm]
\sigma|_{t=s}=\sigma_s(x),\quad v|_{t=s}=v_s(x).
 \end{array}\right.
\end{eqnarray}
By repeating the argument in the proof of Theorem 1.3 in \cite{Y. J. Wang-Z. Tan-2011}, we can get the following result for the problem (\ref{2.6}). The details are omitted here.
\begin{Lemma}
Let $l\geq0$ be an integer. Assume that $(\sigma, v)$ is the
solution of the problem (\ref{2.6}) with the initial date $\sigma_s\in H^{l+1}\cap
L^1$ and $v_s\in H^{l}\cap L^1$, then
\[\|\sigma(t)\|\leq\displaystyle C(1+t)^{-\frac{n}{4}}\left(\|(\sigma_s,v_s)\|_{L^1}+\|(\sigma_s,v_s)\|\right),\]
\[\|\nabla^{k+1}\sigma(t)\|\leq\displaystyle C(1+t)^{-\frac{n}{4}-\frac{k+1}{2}}\left(\|(\sigma_s,v_s)\|_{L^1}+\|(\nabla^{k+1}\sigma_s,\nabla^{k}v_s)\|\right),\]
\[\|\nabla^{k}v(t)\|\leq\displaystyle C(1+t)^{-\frac{n}{4}-\frac{k}{2}}\left(\|(\sigma_s,v_s)\|_{L^1}+\|(\nabla^{k+1}\sigma_s,\nabla^{k}v_s)\|\right),\]
where $k$ is an integer satisfying $0\leq k\leq l$.
\end{Lemma}

\section{Energy estimates}
\setcounter{equation}{0} In this section, we will perform some
energy estimates on solutions $(\sigma, v)$ to problem (\ref{2.3}).
Throughout of this section, we assume that $f(t, x)\in
H^{N-1}(\mathbb{R}^n)\cap L^1(\mathbb{R}^n)$ for all $t\geq0$. For
later use, we list some standard inequalities as follows. cf.
\cite{H. F. Ma-S. Ukai-T. Yang-2010}.
\begin{Lemma}
Let $m$ be a positive integer and $u\in
H^{[\frac{n}{2}]+1}(\mathbb{R}^n)$, then
\[\|u\|^2_{L^\infty}\leq C \|\nabla^{m+1}u\|\|\nabla^{m-1}u\|\quad
for\,\, n=2m,\]
\[\|u\|^2_{L^\infty}\leq C \|\nabla^{m+1}u\|\|\nabla^{m}u\|\quad
for\,\, n=2m+1.\]
\end{Lemma}

\begin{Lemma}
Let $m$ be the integer defined in Lemma 3.1 and $f, g,h\in H^{[\frac{n}{2}]+1}(\mathbb{R}^n)$ , then we have
\[(i) \left|\int_{\mathbb{R}^n}f\cdot g \cdot h\,dx\right|\leq \epsilon\|\nabla^{m-1}f\|_2^2+C_\epsilon\|g\|^2\|h\|^2,\]
\[(ii) \left|\int_{\mathbb{R}^n}f\cdot g \cdot h\,dx\right|\leq \epsilon\|f\|_2^2+C_\epsilon\|\nabla^{m-1}g\|^2_2\|h\|^2,\]
for any $\epsilon>0$. Here and hereafter, $C_\epsilon$ denotes a positive constant depending only on $\epsilon$.
\end{Lemma}
We first give the energy estimate on the low order derivatives of
$(\sigma, v)$.

\begin{Lemma}
Let $n\geq5$, $N\geq n+2$, then there exists two suitably small
constants $d_0>0$ and $\epsilon_0>0$ such that for
$0<\epsilon\leq\epsilon_0$, it holds
\begin{equation}\label{3.1}
\begin{array}{rl}
 &\displaystyle\frac{d}{dt}\left(\|U(t)\|^2+\|\nabla\sigma(t)\|^2+d_0\langle
v, \nabla\sigma\rangle(t)\right)+\|\nabla v(t)\|^2+\|\nabla\sigma(t)\|^2_1\\[2mm]
 &\leq \epsilon C\left(\|\nabla^3\sigma(t)\|^2_{m-2}+\|\nabla^2v(t)\|^2_{m-1}\right)+C_\epsilon C\left(\|\tilde{U}(t)\|^2_{m+1}\|\nabla\tilde{U}(t)\|^2_{1}+\|f(t)\|^2_{L^1\cap
L^2}\right),
\end{array}
\end{equation}
where $m$ is defined in Lemma 3.1 and $C$ depends only on $\rho_\infty, \mu, \nu$ and $\kappa$.
\end{Lemma}
\noindent{\bf Proof.}~~Multiplying $(\ref{2.3})_1$ and
$(\ref{2.3})_2$ by $\sigma$ and $v$, respectively, and integrating
them over $\mathbb{R}^n$, we have from integrating by parts that
\begin{equation}\label{3.2}
\begin{array}{rl}
 &\displaystyle\frac{1}{2}\frac{d}{dt}\|U\|^2+\mu^\prime\|\nabla v\|^2+\nu^\prime\|\nabla\cdot v\|^2\\[2mm]
 =&\langle G_1(\tilde{U}),\sigma\rangle+\langle
 G_2(\tilde{U}),v\rangle+\kappa^\prime\langle\nabla\Delta\sigma,v\rangle+\lambda_2\langle f,v\rangle\\[2mm]
 =&I_0+I_1+I_2+I_3.
\end{array}
\end{equation}
From (\ref{2.2}) and Lemma 3.2, we have
\begin{equation}\label{3.3}
\begin{array}{rl}
 I_0&\leq\epsilon\|\nabla^{m-1}\sigma\|_2^2+C_\epsilon C\left(\|\nabla\tilde{\sigma}\|^2\|\tilde{v}\|^2+\|\tilde{\sigma}\|^2\|\nabla\tilde{v}\|^2\right)\\[2mm]
  &\leq\epsilon\|\nabla^{m-1}\sigma\|_2^2+C_\epsilon C\|\tilde{U}\|^2\|\nabla\tilde{U}\|^2,
\end{array}
\end{equation}
and
\begin{equation}\label{3.4}
 I_1\leq\epsilon\|\nabla^{m-1}v\|^2_2+C_\epsilon C\|\tilde{U}\|^2\|\nabla\tilde{U}\|_1^2.
\end{equation}
For $I_2$, integrating by parts and using $(\ref{2.3})_1$, (\ref{2.2}) and Lemma 3.2, we deduce that
\begin{equation}\label{3.5}
\begin{array}{rl}
 I_2&=\displaystyle-\kappa^\prime\langle\Delta\sigma,\nabla\cdot v\rangle=\frac{\kappa^\prime}{\gamma}\langle\Delta\sigma,\sigma_t-G_1(\tilde{U})\rangle\\[3mm]
 &=\displaystyle-\frac{\kappa^\prime}{2\gamma}\frac{d}{dt}\|\nabla\sigma\|^2-\frac{\kappa^\prime}{\gamma}\langle\Delta\sigma, G_1(\tilde{U})\rangle\\[3mm]
 &\displaystyle\leq-\frac{\kappa^\prime}{2\gamma}\frac{d}{dt}\|\nabla\sigma\|^2+\epsilon\|\nabla^2\sigma\|^2+C_\epsilon C\|\nabla\tilde{U}\|^2\|\nabla^{m-1}\tilde{U}\|^2_2.
 \end{array}
\end{equation}
 For $I_3$, Lemma 3.1 gives
\begin{equation}\label{3.6}
 I_3\leq\epsilon\|\nabla^{m-1}v\|^2_2+C_\epsilon C\|f\|^2_{L^1}.
\end{equation}
Since $n\geq5$, $N\geq n+2$, we have $m-1\geq1$. Substituting
(\ref{3.3})-(\ref{3.6}) into (\ref{3.2}) yields
\begin{equation}\label{3.7}
\begin{array}{rl}
 &\displaystyle\frac{d}{dt}\left(\|U\|^2+\|\nabla\sigma\|^2\right)+\|\nabla v\|^2+\|\nabla\cdot v\|^2\\[2mm]
 &\leq\epsilon C\left(\|\nabla^{m-1}\sigma\|_2^2+\|\nabla^{2}\sigma\|^2\right)+\epsilon
 C\|\nabla^{2}v\|^2_{m-1}+C_\epsilon C\left(\|\tilde{U}\|^2_{m+1}\|\nabla\tilde{U}\|^2_1+\|f\|_{L^1}^2\right),
\end{array}
\end{equation}
provided that $\epsilon$ is small enough, where $C$ depends only on $\rho_\infty, \mu, \nu$ and $\kappa$.

Next, we estimate $\|\nabla\sigma\|^2$. Taking the $L^2$ inner
product with $\nabla\sigma$ on both side of $(\ref{2.3})_2$ and then
integrating by parts, we have
\begin{equation}\label{3.8}
\begin{array}{rl}
 &\displaystyle\gamma\|\nabla\sigma\|^2+\kappa^\prime\|\nabla^2\sigma\|^2\\[2mm]
 &=-\langle v_t,\nabla\sigma\rangle+\mu^\prime\langle
 \Delta v, \nabla\sigma\rangle+\nu^\prime\langle\nabla(\nabla\cdot v), \nabla\sigma\rangle+\langle G_2(\tilde{U})+\lambda_2f, \nabla\sigma\rangle\\[2mm]
 &=I_4+I_5+I_6+I_7.
\end{array}
\end{equation}
Similar to (3.5), the term $I_4$ can be controlled by

\begin{equation}\label{3.9}
\begin{array}{rl}
I_4&=-\displaystyle\frac{d}{dt}\langle v,
\nabla\sigma\rangle-\langle \nabla\cdot v,
\sigma_t\rangle\\[2mm]
&=-\displaystyle\frac{d}{dt}\langle v, \nabla\sigma\rangle-\langle
\nabla\cdot v,
-\gamma\nabla\cdot v+G_1(\tilde{U})\rangle\\[2mm]
&\leq-\displaystyle\frac{d}{dt}\langle v,
\nabla\sigma\rangle+2\gamma\|\nabla\cdot
v\|^2+C\|\nabla^{m-1}\tilde{U}\|^2_2\|\nabla\tilde{U}\|^2.
\end{array}
\end{equation}
Integrating by parts and using the Cauchy-Schwartz inequality, it
is easy to get
\begin{equation}\label{3.10}
 I_5+I_6\leq\frac{\kappa^\prime}{4}\|\nabla^2\sigma\|^2+C(\|\nabla
 v\|^2+\|\nabla\cdot v\|^2).
\end{equation}
Finally, (\ref{2.2}) and the Cauchy-Schwartz inequality imply that
\begin{equation}\label{3.11}
I_7\leq\displaystyle\frac{\gamma}{2}\|\nabla\sigma\|^2+C\left(\|\nabla^{m-1}\tilde{U}\|^2_2\|\nabla
 \tilde{U}\|^2_1+\|f\|^2\right).
\end{equation}
Combining (\ref{3.8})-(\ref{3.11}), we obtain
\begin{equation}\label{3.12}
\begin{array}{rl}
&\displaystyle\frac{d}{dt}\langle v,
\nabla\sigma\rangle+\|\nabla\sigma\|^2+\|\nabla^2\sigma\|^2 \\[2mm]
&\leq C(\|\nabla v\|^2+\|\nabla\cdot
v\|^2)+C\left(\|\nabla^{m-1}\tilde{U}\|^2_2\|\nabla
 \tilde{U}\|^2_1+\|f\|^2\right).
\end{array}
\end{equation}
where the constant $C$ depends only on $\rho_\infty, \mu, \nu$ and $\kappa$. Multiplying (\ref{3.12}) with a small constant $d_0>0$ and then
adding the resultant equation to (\ref{3.7}), one can get
(\ref{3.1}) immediately by the smallness of $d_0$ and $\epsilon$. This completes
the proof of Lemma 3.3.

Next, we derive the energy estimate on the high order derivatives of
$(\sigma, v)$. We establish the following lemma.
\begin{Lemma}
Let $n\geq5$, $N\geq n+2$, then there exists two suitably small
constants $d_1>0$ and $\epsilon_1>0$ such that for
$0<\epsilon\leq\epsilon_1$, it holds
\begin{equation}\label{3.13}
\begin{array}{rl}
 &\displaystyle\frac{d}{dt}\left(\|\nabla\sigma(t)\|^2_N+\|\nabla v(t)\|^2_{N-1}+d_1\sum_{|\alpha|=1}^N\langle
\partial_x^\alpha v, \partial_x^\alpha \nabla\sigma\rangle(t)\right)+\|\nabla^2\sigma(t)\|^2_N+\|\nabla^2 v(t)\|^2_{N-1}\\[4mm]
 &\leq \epsilon C\|\nabla\sigma(t)\|^2+C_\epsilon C\left(\|\nabla\tilde{U}(t)\|^2_{N-2}\|\nabla\tilde{U}(t)\|^2_{N}+\|f(t)\|^2_{N-1}\right),
\end{array}
\end{equation}
where $C$ is depending only on $\rho_\infty, \mu, \nu$ and $\kappa$.
\end{Lemma}
\noindent{\bf Proof.}~~For each multi-index $\alpha$ with
$1\leq|\alpha|\leq N$, applying $\partial_x^\alpha$ to
$(\ref{2.3})_1$ and $(\ref{2.3})_2$ and then taking the $L^2$ inner
product with $\partial_x^\alpha\sigma$ and $\partial_x^\alpha v$ on
the two resultant equations respectively,
 we have from integrating by parts that
 \begin{equation}\label{3.14}
\begin{array}{rl}
 &\displaystyle\frac{1}{2}\frac{d}{dt}\left(\|\partial_x^\alpha\sigma\|^2+\|\partial_x^\alpha v\|^2\right)+\mu^\prime\|\partial_x^\alpha\nabla v\|^2+\nu^\prime\|\partial_x^\alpha\nabla\cdot v\|^2\\[2mm]
 =&\langle \partial_x^\alpha G_1(\tilde{U}),\partial_x^\alpha\sigma\rangle+\langle
 \partial_x^\alpha G_2(\tilde{U}),\partial_x^\alpha v\rangle+\kappa^\prime\langle\partial_x^\alpha\nabla\Delta\sigma,\partial_x^\alpha v\rangle+\lambda_2\langle\partial_x^\alpha f,\partial_x^\alpha v\rangle\\[2mm]
 =&I_8+I_9+I_{10}+I_{11}.
\end{array}
\end{equation}
Now, we estimate $I_8$-$I_{11}$ term by term. For $I_8$, we deduce
from (\ref{2.2}) and the Cauchy-Schwartz inequality that
\begin{equation}\label{3.15}
\begin{array}{rl}
I_8&\leq\epsilon\|\partial_x^\alpha\sigma\|^2+C_\epsilon\|\partial_x^\alpha
 G_1(\tilde{U})\|^2\\[2mm]
 &\leq\epsilon\|\partial_x^\alpha\sigma\|^2+C_\epsilon C
 \left(\|\partial_x^\alpha(\nabla\tilde{\sigma}\cdot\tilde{v})\|^2+\|\partial_x^\alpha(\tilde{\sigma}\nabla\cdot\tilde{v})\|^2\right).
\end{array}
\end{equation}
By Leibniz's formula and Minkowski's inequality, we get
\begin{equation}\label{3.16}
\begin{array}{rl}
\|\partial_x^\alpha(\nabla\tilde{\sigma}\cdot\tilde{v})\|^2\leq&
C(\|(\partial_x^\alpha\nabla\tilde{\sigma})\cdot\tilde{v}\|^2+\|\nabla\tilde{\sigma}\cdot\partial_x^\alpha\tilde{v}\|^2)
+C\displaystyle\sum_{0<|\beta|=|\alpha|-1}C^\alpha_\beta\|\partial_x^\beta\nabla\tilde{\sigma}\cdot\partial_x^{\alpha-\beta}\tilde{v}\|^2\\[2mm]
&+C\displaystyle\sum_{0<|\beta|\leq|\alpha|-2,\,
|\alpha-\beta|\leq\frac{N}{2}}C^\alpha_\beta\|\partial_x^\beta\nabla\tilde{\sigma}\cdot\partial_x^{\alpha-\beta}\tilde{v}\|^2\\[2mm]
&+C\displaystyle\sum_{0<|\beta|\leq|\alpha|-2,\,
|\alpha-\beta|>\frac{N}{2}}C^\alpha_\beta\|\partial_x^\beta\nabla\tilde{\sigma}\cdot\partial_x^{\alpha-\beta}\tilde{v}\|^2\\[2mm]
=&J_0+J_1+J_2+J_3.
\end{array}
\end{equation}
Here $C^\alpha_\beta$ denotes the binomial coefficients
corresponding to multi-indices.  For $J_0$, lemma 3.1 gives
\begin{equation}\label{3.17}
\begin{array}{rl}
J_0\leq&
C\left(\|\tilde{v}\|^2_{L^\infty}\|\partial_x^\alpha\nabla\tilde{\sigma}\|^2+\|\nabla\tilde{\sigma}\|^2_{L^\infty}\|\partial_x^\alpha\tilde{v}\|^2\right)\\[2mm]
\leq&
C\left(\|\nabla\tilde{v}\|^2_{N-5}\|\nabla^2\tilde{\sigma}\|^2_{N-1}+\|\nabla^2\tilde{\sigma}\|^2_{N-5}\|\nabla\tilde{v}\|^2_{N-1}\right),
\end{array}
\end{equation}
where, in the last inequality of (\ref{3.17}), we have used the fact
that  $m-1\geq1$ and $m+1\leq N-4$ due to $N\geq n+2$ and $n\geq5$.
Similarly, it holds that
\begin{equation}\label{3.18}
J_1\leq
C\displaystyle\sum_{0<|\beta|=|\alpha|-1}\|\partial_x^{\alpha-\beta}\tilde{v}\|^2_{L^\infty}\|\partial_x^\beta\nabla\tilde{\sigma}\|^2\\[3mm]
\leq
C\|\nabla^2\tilde{v}\|^2_{N-5}\|\nabla^2\tilde{\sigma}\|_{N-2}^2.
\end{equation}
For the terms $J_2$ and $J_3$, notice that for any $\beta\leq\alpha$
with $|\alpha-\beta|\leq\frac{N}{2}$,
\[|\alpha-\beta|+m+1\leq\frac{N}{2}+\frac{n}{2}+1\leq\frac{N}{2}+\frac{N}{2}=N,\]
and for any $\beta\leq\alpha$ with $|\alpha-\beta|>\frac{N}{2}$,
\[|\beta|+m+2=|\alpha|-|\alpha-\beta|+m+2<N-\frac{N}{2}+\frac{n}{2}+2\leq N+1.\]
which implies  $|\beta|+m+2\leq N$ since $|\beta|$ and $m$ are
positive integers. Hence, we deduce from Lemma 3.1 that
\begin{equation}\label{3.19}
J_2\leq C\displaystyle\sum_{0<|\beta|\leq|\alpha|-2,\,
|\alpha-\beta|\leq\frac{N}{2}}\|\partial_x^{\alpha-\beta}\tilde{v}\|^2_{L^\infty}\|\partial_x^\beta\nabla\tilde{\sigma}\|^2
\leq
C\|\nabla^2\tilde{v}\|^2_{N-2}\|\nabla^2\tilde{\sigma}\|_{N-3}^2,
\end{equation}
and
\begin{equation}\label{3.20}
J_3\leq C\displaystyle\sum_{0<|\beta|\leq|\alpha|-2,\,
|\alpha-\beta|>\frac{N}{2}}\|\partial_x^\beta\nabla\tilde{\sigma}\|^2_{L^\infty}\|\partial_x^{\alpha-\beta}\tilde{v}\|^2
\leq
C\|\nabla^2\tilde{v}\|^2_{N-3}\|\nabla^2\tilde{\sigma}\|_{N-2}^2.
\end{equation}
Putting (\ref{3.17})-(\ref{3.20}) into (\ref{3.16}), we arrive at
\begin{equation}\label{3.21}
\|\partial_x^\alpha(\nabla\tilde{\sigma}\cdot\tilde{v})\|^2\leq
C\left(\|\nabla\tilde{v}\|^2_{N-5}\|\nabla^2\tilde{\sigma}\|^2_{N-1}+\|\nabla^2\tilde{U}\|^2_{N-3}\|\nabla\tilde{U}\|^2_{N-1}\right).
\end{equation}
Similarly, it holds
\begin{equation}\label{3.22}
\|\partial_x^\alpha(\tilde{\sigma}\nabla\cdot\tilde{v})\|^2\leq
C\left(\|\nabla\tilde{\sigma}\|^2_{N-5}\|\nabla^2\tilde{v}\|^2_{N-1}+\|\nabla^2\tilde{U}\|^2_{N-3}\|\nabla\tilde{U}\|^2_{N-1}\right).
\end{equation}
Combining (\ref{3.15}), (\ref{3.21}) and (\ref{3.22}) yields
\begin{equation}\label{3.23}
I_8\leq\epsilon\|\partial_x^\alpha\sigma\|^2+C_\epsilon C
 \|\nabla\tilde{U}\|^2_{N-2}\|\nabla\tilde{U}\|^2_{N}.
\end{equation}
For the term $I_9$, let $\alpha_0\leq\alpha$ with $|\alpha_0|=1$,
then
\begin{equation}\label{3.24}
I_9=-\langle\partial_x^{\alpha-\alpha_0}G_2,\partial_x^{\alpha+\alpha_0}v\rangle\leq\epsilon\|\partial_x^{\alpha+\alpha_0}v\|^2+C_\epsilon\|\partial_x^{\alpha-\alpha_0}G_2\|^2.
\end{equation}
Similar to the estimate of (\ref{3.21}), we have
\begin{equation}\label{3.25}
\|\partial_x^{\alpha-\alpha_0}G_2\|^2\leq
C\|\tilde{U}\|^2_{N-1}\|\nabla\tilde{ U}\|^2_N.
\end{equation}
Thus, it follows from (\ref{3.24}) and (\ref{3.25}) that
\begin{equation}\label{3.26}
I_9\leq
\epsilon\|\partial_x^{\alpha+\alpha_0}v\|^2+C\|\tilde{U}\|^2_{N-1}\|\nabla\tilde{
U}\|^2_N.
\end{equation}
Notice that (\ref{3.21}) and (\ref{3.22}) imply
\begin{equation}\label{3.27}
\|\partial_x^{\alpha}G_1\|^2\leq C\left(\|\partial_x^\alpha(\nabla\tilde{\sigma}\cdot\tilde{v})\|^2+\|\partial_x^\alpha(\tilde{\sigma}\nabla\cdot\tilde{v})\|^2\right)\leq
C\|\nabla\tilde{U}\|^2_{N-2}\|\nabla\tilde{ U}\|^2_N.
\end{equation}
Therefore, we derive from $(\ref{2.3})_1$, (\ref{3.27}) and the Cauchy-Schwartz inequality that
\begin{equation}\label{3.28}
\begin{array}{rl}
I_{10}=&-\displaystyle\frac{\kappa^\prime}{\gamma}\langle\partial_x^\alpha\Delta\sigma,
-\partial_x^\alpha\sigma_t+\partial_x^\alpha
G_1(\tilde{U})\rangle\\[3mm]
=&-\displaystyle\frac{\kappa^\prime}{\gamma}\langle\partial_x^\alpha\nabla\sigma,
\partial_x^\alpha\nabla\sigma_t\rangle-\displaystyle\frac{\kappa^\prime}{\gamma}\langle\partial_x^\alpha\Delta\sigma,\partial_x^\alpha
G_1(\tilde{U})\rangle\\[2mm]
\leq&-\displaystyle\frac{\kappa^\prime}{2\gamma}\displaystyle\frac{d}{dt}\|\partial_x^\alpha\nabla\sigma\|^2
+\epsilon\|\partial_x^\alpha\Delta\sigma\|^2+C_\epsilon C\|\nabla\tilde{U}\|^2_{N-2}\|\nabla\tilde{U}\|^2_{N}.
\end{array}
\end{equation}
Moreover, it holds that
\begin{equation}\label{3.29}
I_{11}=-\lambda_2\langle\partial_x^{\alpha+\alpha_0}v,\partial_x^{\alpha-\alpha_0}f\rangle\leq\epsilon\|\partial_x^{\alpha+\alpha_0}v\|^2+C_\epsilon C\|f\|^2_{N-1}.
\end{equation}
where $\alpha_0$ is defined in (\ref{3.24}). Combining (\ref{3.14}),
(\ref{3.23}), (\ref{3.26}), (\ref{3.28}) and (\ref{3.29}), if
$\epsilon$ is small enough, we have
\begin{equation}\label{3.30}
\begin{array}{rl}
 &\displaystyle\frac{d}{dt}\left(\|\partial_x^\alpha\sigma\|^2_1+\|\partial_x^\alpha v\|^2\right)+\|\partial_x^\alpha\nabla v\|^2+\|\partial_x^\alpha\nabla\cdot v\|^2\\[2mm]
 &\leq\epsilon C\|\partial_x^\alpha\sigma\|^2+\epsilon
 C\|\partial_x^\alpha\Delta \sigma\|^2+C_\epsilon C\left(\|\nabla\tilde{U}\|^2_{N-2}\|\nabla\tilde{U}\|^2_N+\|f\|_{N-1}^2\right),
\end{array}
\end{equation}
where $C$ depends only on $\rho_\infty, \mu, \nu$ and $\kappa$.

Now we turn to estimate $\|\partial_x^\alpha\Delta\sigma\|^2$ for
$1\leq|\alpha|\leq N$. As we did for the first order derivative
estimate, applying $\partial_x^\alpha$ to $(\ref{2.3})_2$ and then
taking the $L^2$ inner product with $\partial_x^\alpha\nabla\sigma$
on the resultant equation, we get from integrating by parts that
\begin{equation}\label{3.31}
\begin{array}{rl}
 &\displaystyle\kappa^\prime\|\partial_x^\alpha\Delta\sigma\|^2+\gamma\|\partial_x^\alpha\nabla\sigma\|^2\\[2mm]
 &=-\langle \partial_x^\alpha v_t, \partial_x^\alpha\nabla\sigma\rangle+\mu^\prime\langle
 \partial_x^\alpha\Delta v, \partial_x^\alpha\nabla\sigma\rangle+
 \nu^\prime\langle\partial_x^\alpha\nabla(\nabla\cdot v), \partial_x^\alpha\nabla\sigma\rangle\\[2mm]
  &\hspace{4mm}+\langle \partial_x^\alpha
 G_2(\tilde{U}),\partial_x^\alpha\nabla\sigma\rangle
 +\lambda_2\langle \partial_x^\alpha f, \partial_x^\alpha\nabla\sigma\rangle\\[2mm]
 &=I_{12}+I_{13}+I_{14}+I_{15}+I_{16}.
\end{array}
\end{equation}
The first term $I_{12}$ is controlled by
\begin{equation}\label{3.32}
\begin{array}{rl}
I_{12}=&-\displaystyle\frac{d}{dt}\langle\partial_x^\alpha v,
\partial_x^\alpha\nabla\sigma\rangle
+\langle\partial_x^\alpha v, \partial_x^\alpha\nabla\sigma_t\rangle\\[3mm]
=&-\displaystyle\frac{d}{dt}\langle\partial_x^\alpha v,
\partial_x^\alpha\nabla\sigma\rangle-\displaystyle\langle\partial_x^\alpha \nabla\cdot
v, \partial_x^\alpha(-\gamma\nabla\cdot v+G_1(\tilde{U}))\rangle
\\[2mm]
\leq&-\displaystyle\frac{d}{dt}\langle\partial_x^\alpha v,
\partial_x^\alpha\nabla\sigma\rangle+2\gamma\|\partial_x^\alpha
\nabla\cdot v \|^2+C\|\nabla\tilde{U}\|^2_{N-2}\|\nabla\tilde{
U}\|^2_N.
\end{array}
\end{equation}
Here, in the last inequality of (\ref{3.32}), we have used
(\ref{3.27}). By integrating by parts, the Cauchy-Schwartz
inequality and (\ref{3.25}), the other terms $I_{13}$-$I_{15}$ can be
estimated as follows.
\begin{equation}\label{3.33}
 I_{13}+I_{14}\leq\frac{\kappa^\prime}{4}\|\partial_x^\alpha\nabla^2\sigma\|^2+C\left(\|\partial_x^\alpha\nabla
 v\|^2+\|\partial_x^\alpha\nabla\cdot v\|^2\right),
\end{equation}
\begin{equation}\label{3.34}
 I_{15}\leq\frac{\kappa^\prime}{4}\|\partial_x^{\alpha+\alpha_0}\nabla\sigma\|^2+C\|\tilde{U}\|^2_{N-1}\|\nabla\tilde{
 U}\|^2_N,
\end{equation}
\begin{equation}\label{3.35}
 I_{16}\leq\frac{\kappa^\prime}{4}\|\partial_x^{\alpha+\alpha_0}\nabla\sigma\|^2+C\|f\|^2_{N-1}.
\end{equation}
where $\alpha_0$ is given in (\ref{3.24}). Combining
(\ref{3.31})-(\ref{3.35}), we obtain
\begin{equation}\label{3.36}
\begin{array}{rl}
 &\displaystyle\frac{d}{dt}\langle\partial_x^\alpha v,
\partial_x^\alpha\nabla\sigma\rangle+\kappa^\prime\|\partial_x^\alpha\Delta\sigma\|^2+\gamma\|\partial_x^\alpha\nabla\sigma\|^2\\[2mm]
 &\leq C\left(\|\partial_x^\alpha \nabla
v\|^2+\|\partial_x^\alpha \nabla\cdot
v\|^2\right)+C\left(\|\tilde{U}\|^2_{N-1}\|\nabla\tilde{U}\|^2_N+\|f\|_{N-1}^2\right).
\end{array}
\end{equation}
Multiplying (\ref{3.36}) with a suitably small constant $d_1>0$ and then
adding the resultant equation to (\ref{3.30}) gives
\begin{equation}\label{3.37}
\begin{array}{rl}
 &\displaystyle\frac{d}{dt}\left(\|\partial_x^\alpha\sigma\|^2_1+\|\partial_x^\alpha v\|^2+d_1\langle\partial_x^\alpha v,
\partial_x^\alpha\nabla\sigma\rangle\right)+\|\partial_x^\alpha\nabla\sigma\|^2_1+\|\partial_x^\alpha\nabla v\|^2\\[2mm]
 &\leq \epsilon C\|\partial_x^\alpha\sigma\|^2+CC_\epsilon\left(\|\nabla\tilde{U}\|^2_{N-2}\|\nabla\tilde{U}\|^2_N+\|f\|_{N-1}^2\right),
\end{array}
\end{equation}
provided that $d_1$ and $\epsilon$ are small enough, where $C$ depends only on $\rho_\infty, \mu, \nu$ and $\kappa$. Summing up
$\alpha$ with $1\leq|\alpha|\leq N$ in (\ref{3.37}), then
(\ref{3.13}) follows immediately by the smallness of $\epsilon$. This completes the
proof of Lemma 3.4.

As a consequence of Lemmas 3.3-3.4, we have the following Corollary.
\begin{Corollary}
Let $n\geq5$, $N\geq n+2$, then there exists two suitably small
constants  $d_0>0$  and $d_1>0$ such that
\begin{equation}\label{3.38}
\begin{array}{rl}
 &\displaystyle\frac{d}{dt}\left(\|\sigma(t)\|^2_{N+1}+\|v(t)\|^2_{N}+d_0\langle
v,\nabla\sigma\rangle(t)+d_1\displaystyle\sum_{|\alpha|=1}^{N}\langle\partial_x^\alpha
v, \partial_x^\alpha\nabla\sigma\rangle(t)\right)+\|\nabla\sigma(t)\|^2_{N+1}+\|\nabla v(t)\|^2_{N}\\[4mm]
 &\leq
 C\left(\|\tilde{U}(t)\|^2_{N-1}\|\nabla\tilde{U}(t)\|^2_N+\|f(t)\|^2_{H^{N-1}\cap{L^1}}\right),
\end{array}
\end{equation}
where $C$ depends only on $\rho_\infty, \mu, \nu$ and $\kappa$.
\end{Corollary}
\noindent{\bf Proof.}~~Notice that, from the fact that $m-1\geq1$ and $m+1\leq N-4$, we have
$$\|\nabla^3\sigma\|^2_{m-2}+\|\nabla^2v\|^2_{m-1}\leq
C\|\nabla^2\tilde{U}\|_{N-6}^2,$$ and
$$\|\tilde{U}\|_{m+1}\leq\|\tilde{U}\|_{N-4}^2.$$ Adding (\ref{3.37})
to (\ref{3.1}), we obtain (\ref{3.38}) immediately by the smallness of
$\epsilon$. This completes the proof of Corollary 3.1.

\section{Existence of time periodic solution}
\setcounter{equation}{0}
In this section, we will combine the linearized decay estimate Lemma 2.1 with the energy estimates
Corollary 3.1 to show the existence of time periodic solution to (1.1). Now, we are ready to prove Theorem 1.1 as
follows.\\
 \noindent{\bf Proof of Theorem 1.1.}~~The proof is divided into
 two steps.

 {\it Step 1.}~~Suppose that there exists a time periodic
 solution $U^{per}(t):=(\sigma^{per}(x,t), v^{per}(x,t)), t\in\mathbb{R}$ of the system (\ref{2.1}) with
 period $T$, and $U^{per}(t)\in X_{M_0}(0,T)$ for some constant $M_0>0$.  Then it solves (\ref{2.2a}) with initial
 date $U_s=U^{per}(s)$ for any given time $s\in\mathbb{R}$. Choosing
 $s=-kT$ for $k\in\mathbb{N}$. Clearly, $U^{per}(-kT)=U^{per}(0)$, thus (\ref{2.2a})
 can be written in the mild form as
\begin{equation}\label{4.1}
U^{per}(t)=\displaystyle\mathbb{S}(t,
-kT)U^{per}(0)+\int_{-kT}^t\mathbb{S}(t,\tau)(G(U^{per})(\tau)+F(\tau))d\tau.
\end{equation}
Denote $\mathbb{S}(t, -kT)U^{per}(0):=(\sigma^{per}_1(t), v^{per}_1(t))$. Applying Lemma 2.1 to $\mathbb{S}(t, -kT)U^{per}(0)$, we have
\begin{equation}\label{4.2}
\begin{array}{rl}
\|\sigma^{per}_1(t)\|_N\leq&\displaystyle(1+t+kT)^{-\frac{n}{4}}\left(\|(\sigma_0^{per},
v_0^{per})\|_{L^1}+\|\sigma_0^{per}\|_N^2+\|v_0^{per}\|^2_{N-1}\right)\\[2mm]
&\longrightarrow0\quad as \quad k\rightarrow\infty.
\end{array}
\end{equation}
and
\begin{equation}\label{4.3}
\begin{array}{rl}
\|v_1^{per}(t)\|_{N-1}\leq&\displaystyle(1+t+kT)^{-\frac{n}{4}}\left(\|(\sigma_0^{per},
v_0^{per})\|_{L^1}+\|\sigma_0^{per}\|_N^2+\|v_0^{per}\|^2_{N-1}\right)\\[2mm]
&\longrightarrow0\quad as \quad k\rightarrow\infty.
\end{array}
\end{equation}
Since $L^2\cap L^1$ is dense in $L^2$, (\ref{4.2}) and (\ref{4.3})
still hold for $U^{per}(0)=(\sigma^{per}_0, v^{per}_0)\in
H^N(\mathbb{R}^n)\times H^{N-1}(\mathbb{R}^n)$.  On the other hand,
denote
\[\mathbb{S}(t,\tau)(G(U^{per})(\tau)+F(\tau)):=(S_1(t,\tau),
S_2(t,\tau)).\]
By using Lemma 2.1 again, we get
\begin{equation}\label{4.4}
\|S_1(t,\tau)\|_N\leq\displaystyle(1+t-\tau)^{-\frac{n}{4}}K_0,\quad\|S_2(t,\tau)\|_{N-1}\leq\displaystyle(1+t-\tau)^{-\frac{n}{4}}K_0,
\end{equation}
where \[
\begin{array}{rl}
 K_0=&\|(G_1(U^{per}),
G_2(U^{per})+\lambda_2f)(\tau)\|_{L^1}\\[2mm]
&+\|G_1(U^{per})(\tau)\|_{N}+\|(G_2(U^{per})+\lambda_2f)(\tau)\|_{N-1}.
\end{array}
\]
Then (\ref{4.4}) guarantees the convergence of the integral in
(\ref{4.1}) since $\frac{n}{4}>1$ when $n\geq5$. Thus, letting
$k\rightarrow\infty$ in (\ref{4.1}), we obtain
\begin{equation}\label{4.5}
U^{per}(t)=\displaystyle\int_{-\infty}^t\mathbb{S}(t,\tau)(G(U^{per})+F)(\tau)d\tau.
\end{equation}
For any $U=(\sigma, v)\in X_{M_0}(0,T)$, define
\[\Psi[U](t)=\displaystyle\int_{-\infty}^t\mathbb{S}(t,\tau)(G(U)+F)(\tau)d\tau.\]
Then (\ref{4.5}) shows that $U^{per}$ is a fixed point of $\Psi[U]$.

Conversely, suppose that $\Psi$ has a unique fixed point, denoted by
$U_1(t)=(\sigma_1, v_1)(t)$. We show that $U_1(t)$ is time periodic with period $T$. To this end,
setting $U_2(t)=U_1(t+T)$.  Since the period of $f$ is $T$, the period of $F$ is $T$ too.
Thus, we have
\begin{equation}\label{4.6}
\begin{array}{rl}
U_2(t)&=U_1(t+T)=\Psi[U_1](t+T)\\[3mm]
      &=\displaystyle\int^{t+T}_{-\infty}\mathbb{S}(t+T,\tau)(G(U_1)(\tau)+F(\tau))d\tau\\[3mm]
      &=\displaystyle\int^{t}_{-\infty}\mathbb{S}(t+T, s+T)\left(G(U_1)(s+T)+F(s+T)\right)ds\\[3mm]
      &=\displaystyle\int^{t}_{-\infty}\mathbb{S}(t, s)\left(G(U_2)(s)+F(s)\right)ds\\[3mm]
      &=\Psi[U_2](t)
\end{array}
\end{equation}
where we have used
$$\mathbb{S}(t+T, s+T)=\mathbb{S}(t, s).$$
Then by uniqueness, $U_2=U_1$, which proves the periodicity of
$U_1(t)$. Since $U_1(t)$ is differentiable with respect to $t$, it is the desired periodic solution of the system (2.1).

{\it Step 2.}~~Now, it remains to show that if (H1)-(H3) hold, and
\[ \sup_{0\leq t\leq T}\|f(t)\|_{H^{N-1}\cap L^1}\] is sufficiently
small, then $\Psi$ has a unique fixed point in the space
$X_{M_0}(0,T)$ for some appropriate constant $M_0>0$. The proof is
divided into two parts.

(i)~Assume that $\tilde{U}=(\tilde{\sigma},
 \tilde{v})$ in the system (\ref{2.3}) is time periodic with period $T$. Denote $U=\Psi[\tilde{U}]$ with $U=(\sigma, v)$. Then by the same argument as (\ref{4.6}), one can show that $U$ is also time
 periodic with period $T$.  Notice that $U$ satisfies
 the system (\ref{2.3}). Thus, for $n\geq5$ and $N\geq n+2$, Corollary 3.1 holds.
 Integrating (\ref{3.38}) in $t$ over $[0,T]$ to get
\begin{equation}\label{4.7}
\begin{array}{rl}
&\displaystyle\int^{T}_{0}\left(\|\nabla\sigma(t)\|^2_{N+1}+\|\nabla v(t)\|^2_{N}\right)dt\\[3mm]
      &\leq C\displaystyle\int^{T}_{0}\left(\|\tilde{U}(t)\|^2_{N-1}\|\nabla\tilde{U}(t)\|^2_N+\|f(t)\|_{N-1}^2+\|f(t)\|^2_{L^1}\right)dt\\[3mm]
&\leq C\displaystyle\sup_{0\leq t\leq
T}\|\tilde{U}(t)\|^2_{N-1}\displaystyle\int^{T}_{0}\|\nabla\tilde{U}(t)\|^2_Ndt+\displaystyle\int^{T}_{0}\|f(t)\|^2_{H^{N-1}\cap{L^1}}dt\\[3mm]
&\leq C |||\tilde{U}(t)|||^4+CT\displaystyle\sup_{0\leq t\leq
T}\|f(t)\|_{H^{N-1}\cap{L^1}}^2.
\end{array}
\end{equation}
On the other hand, by Lemma 2.1, we have
\begin{equation}\label{4.8}
\|\sigma(t)\|_N\leq\displaystyle\int^t
_{-\infty}(1+t-\tau)^{-\frac{n}{4}}K_1\,d\tau,\quad\|v(t)\|_{N-1}\leq\displaystyle\int^t
_{-\infty}(1+t-\tau)^{-\frac{n}{4}}K_1\,d\tau,
\end{equation}
where\begin{equation}\label{4.9}
\begin{array}{rl}
 K_1=&\|(G_1(\tilde{U}),
G_2(\tilde{U})+\lambda_2f)(\tau)\|_{L^1}\\[2mm]
&+\|G_1(\tilde{U})(\tau)\|_{N}+\|(G_2(\tilde{U})+\lambda_2f)(\tau)\|_{N-1}.
\end{array}
\end{equation}
From (\ref{2.2}), (\ref{3.25}) and (\ref{3.27}), we easily deduce that
\begin{equation}\label{4.10}
\begin{array}{rl}
\|(G_1(\tilde{U})(\tau)\|_{L^1}&\leq
C\|\nabla\tilde{U}(\tau)\|\|\tilde{U}(\tau)\|,\\[2mm]
\|(G_1(\tilde{U})(\tau)\|_N&\leq
C\|\nabla\tilde{U}(\tau)\|_{N-2}\|\nabla\tilde{U}(\tau)\|_N,\\[2mm]
\|G_2(\tilde{U})+\lambda_2f)(\tau)\|_{L^1}&\leq C\|\nabla\tilde{U}(\tau)\|_1\|\tilde{U}(\tau)\|+C\|f(\tau)\|_{L^1},\\[2mm]
\|G_2(\tilde{U})+\lambda_2f)(\tau)\|_{N-1}&\leq
C\|\tilde{U}(\tau)\|_{N-1}\|\nabla\tilde{U}(\tau)\|_N+C\|f(\tau)\|_{N-1}.
\end{array}
\end{equation}
Combining (\ref{4.8})-(\ref{4.10}), we obtain
\begin{equation}\label{4.11}
\begin{array}{rl}
\|\sigma(t)\|_N\leq&
C\displaystyle\int^t_{-\infty}(1+t-\tau)^{-\frac{n}{4}}\left(\|\tilde{U}(\tau)\|_{N-1}\|\nabla\tilde{U}(\tau)\|_N+\|f(\tau)\|_{H^{N-1}\cap{L^1}}\right)d\tau\\[2mm]
\leq&
C\displaystyle\sum_{j=0}^{\infty}A_j+C\int^t_{-\infty}(1+t-\tau)^{-\frac{n}{4}}\|f(\tau)\|_{H^{N-1}\cap{L^1}}d\tau\\[2mm]
\leq&C\displaystyle\sum_{j=0}^{\infty}A_j+C \sup_{0\leq t\leq
T}\|f(t)\|_{H^{N-1}\cap{L^1}},
\end{array}
\end{equation}
where
\begin{equation}\label{4.12}
\begin{array}{rl}
A_j=&
C\displaystyle\int^{t-jT}_{t-(j+1)T}(1+t-\tau)^{-\frac{n}{4}}\|\tilde{U}(\tau)\|_{N-1}\|\nabla\tilde{U}(\tau)\|_Nd\tau\\[2mm]
\leq&
C\displaystyle\left(\int^{t-jT}_{t-(j+1)T}(1+t-\tau)^{-\frac{n}{2}}d\tau\right)^{\frac{1}{2}}\left(\int^{t-jT}_{t-(j+1)T}\|\tilde{U}(\tau)\|_{N-1}^2\|\nabla\tilde{U}(\tau)\|_N^2d\tau\right)^{\frac{1}{2}}\\[2mm]
\leq&C \displaystyle(1+jT)^{-\frac{n}{4}}\sup_{0\leq \tau\leq
T}\|\tilde{U}(\tau)\|_{N-1}\left(\int^{T}_{0}\|\nabla\tilde{U}(\tau)\|_N^2d\tau\right)^{\frac{1}{2}}\\[3mm]
\leq&C \displaystyle(1+jT)^{-\frac{n}{4}}|||\tilde{U}|||^2
\end{array}
\end{equation}
Since $\frac{n}{4}>1$ when $n\geq5$, substituting (\ref{4.12}) into
(\ref{4.11}) gives
\begin{equation}\label{4.13}
\|\sigma(t)\|_N\leq C |||\tilde{U}|||^2+C \sup_{0\leq t\leq
T}\|f(t)\|_{H^{N-1}\cap{L^1}}.
\end{equation}
Similarly, it holds that
\begin{equation}\label{4.14}
\|v(t)\|_{N-1}\leq C |||\tilde{U}|||^2+C \sup_{0\leq t\leq
T}\|f(t)\|_{H^{N-1}\cap{L^1}}.
\end{equation}
Thus, we deduce from (\ref{4.7}), (\ref{4.13}) and (\ref{4.14}) that
\begin{equation}\label{4.15}
|||\Psi[\tilde{U}]|||\leq C_1 |||\tilde{U}|||^2+C_2 \sup_{0\leq
t\leq T}\|f(t)\|_{H^{N-1}\cap{L^1}},
\end{equation}
where $C_1$ and $C_2$ are some positive constants depending only on $\rho_\infty, \mu, \nu, \kappa$ and $T$.

(ii)~Let $\tilde{U}_1=(\tilde{\sigma}_1, \tilde{v}_1)$ and
$\tilde{U}_2=(\tilde{\sigma}_2, \tilde{v}_2)$ be time periodic functions with period $T$ in the
space $X_{M_0}(0,T)$, where $M_0>0$ will be determined below.  Then
similar to (i), we can get
\begin{equation}\label{4.16}
|||\Psi[\tilde{U}_1]-\Psi[\tilde{U}_2]|||\leq C_3
\left(|||\tilde{U}_1|||+|||\tilde{U}_2|||\right)|||\tilde{U}_1-\tilde{U}_2|||,
\end{equation}
where $C_3$ is a positive constant depending only on $\rho_\infty, \mu, \nu, \kappa$ and $T$. Choose $M_0>0$ and a sufficiently small constant
$\delta>0$  such that
\begin{equation}\label{4.17}
 C_1M_0^2+C_2\delta\leq M_0,\quad\mbox{and}\,\,2C_3M_0<1
\end{equation}
That is,
\begin{equation}\label{4.18}
\frac{1-\sqrt{1-4C_1C_2\delta}}{2C_1}\leq
M_0\leq\min\left\{\frac{1+\sqrt{1-4C_1C_2\delta}}{2C_1},
\frac{1}{2C_3}, 1\right\}
\end{equation}
Notice that
\[
\frac{1-\sqrt{1-4C_1C_2\delta}}{2C_1}\longrightarrow0\quad as\quad
\delta\longrightarrow0.
\]
Then there exists a constant $\delta_0>0$ depending only on $\rho_\infty, \mu, \nu, \kappa$ and $T$ such that if
$0<\delta\leq\delta_0$, the set of $M_0$ that satisfying
(\ref{4.18}) is not empty. For $0<\delta\leq\delta_0$, when $M_0$
satisfies (\ref{4.18}), $\Psi$ is a contraction map in the complete
space $X_{M_0}(0,T)$, thus $\Psi$ has a unique fixed point in
$X_{M_0}(0,T)$. This completes the proof of Theorem 1.1.

\section{Stability of time periodic solution}
\setcounter{equation}{0} This section is devoted to proving Theorem
1.2 on the stability of the obtained time periodic solution. We
shall establish the global existence of smooth solutions to the
Cauchy problem (1.1), (1.5).

First, let $(\rho^{per}, u^{per})$ be the time periodic solution
obtained in Theorem 1.1 and $(\rho, u)$ be the solution of the
Cauchy problem (1.1), (1.5). Denote
\[
(\sigma^{per}, v^{per})=(\rho^{per}-\rho_\infty, \lambda_2u^{per}),
\]
\[
(\sigma, v)=(\rho-\rho_\infty, \lambda_2u).
\]
Let $(\bar{\sigma}, \bar{v})=(\sigma-\sigma^{per}, v-v^{per})$, then
$(\bar{\sigma}, \bar{v})$ satisfies
 \begin{eqnarray}\label{5.1}
\left\{\begin{array}{ll}
          \bar{\sigma}_t+\gamma\nabla\cdot \bar{v}=G_1(\bar{\sigma}+\sigma^{per}, \bar{v}+v^{per})-G_1(\sigma^{per}, v^{per}), \\[2mm]
         \bar{ v}_t-\mu^\prime\Delta \bar{v}-\nu^\prime\nabla(\nabla\cdot
          \bar{v})+\gamma\nabla\bar{\sigma}-\kappa^\prime\nabla\Delta\bar{\sigma}=G_2(\bar{\sigma}+\sigma^{per}, \bar{v}+v^{per})-G_2(\sigma^{per}, v^{per}),
 \end{array}\right.
\end{eqnarray}
with the initial date
 \begin{eqnarray}\label{5.2}
\bar{\sigma}|_{t=0}=\bar{\sigma}_0(x)=\rho_0(x)-\rho^{per}(0),\quad
\bar{v}|_{t=0}=\bar{v}_0(x)=\lambda_2(u_0(x)-u^{per}(0)).
\end{eqnarray}
Define the solution space by $\bar{X}(0, \infty)$, where for $0\leq
t_1\leq t_2\leq\infty$,
\begin{equation}\label{5.3}
\bar{X}(t_1, t_2)=\left\{(\bar{\sigma}, \bar{v})(t,x)\left|
\begin{array}{c}
\bar{\sigma}(t,x)\in C(t_1, t_2; H^{N-1}(\mathbb{R}^n))\cap C^1(t_1, t_2; H^{N-3}(\mathbb{R}^n)),\\[2mm]
\bar{v}(t,x)\in C(t_1, t_2; H^{N-2}(\mathbb{R}^n))\cap C^1(t_1, t_2; H^{N-4}(\mathbb{R}^n)),\\[2mm]
\nabla\bar{\sigma}(t,x)\in L^2(t_1, t_2; H^{N-1}(\mathbb{R}^n)),
\nabla\bar{v}(t,x)\in L^2(t_1, t_2; H^{N-2}(\mathbb{R}^n)),
\end{array}
\right.\right\}
\end{equation}
with the norm
\begin{equation}\label{5.4}
\|(\bar{\sigma}, \bar{v})(t)\nparallel^2:=\sup_{t_1\leq t\leq
t_2}\left\{\|\bar{\sigma}(t)\|_{N-1}^2+\|\bar{v}(t)\|_{N-2}^2\right\}+\int_{t_1}^{t_2}\left(\|\nabla\bar{\sigma}(t)\|_{N-1}^2+\|\nabla
\bar{v}(t)\|_{N-2}^2\right)dt.
\end{equation}
Notice that $(\sigma^{per},v^{per})\in \bar{X}(0,T)$.

By using the dual argument and iteration technique as \cite{H.
Hattori-D. Li-1994}, one can prove the following local existence of
the Cauchy problem (\ref{5.1}), (\ref{5.2}). We omit the proof
here for brevity.
\begin{Lemma}(Local existence)
Under the assumptions of Theorem 1.1, suppose that
$(\bar{\sigma}_0,\bar{v}_0)\in H^{N-1}(\mathbb{R}^n)\times
H^{N-2}(\mathbb{R}^n)$ and $\inf\rho_0(x)>0$. Then there exists a
positive constant $T_0$ depending only on $\|(\bar{\sigma}_0,
\bar{v}_0)\nparallel$ such that the Cauchy problem (\ref{5.1}),
(\ref{5.2}) admits a unique classical solution $(\bar{\sigma},\bar{v})\in
\bar{X}(0,T_0)$ which satisfies
\[\|(\bar{\sigma},
\bar{v})(t)\nparallel\leq C_4\|(\bar{\sigma}_0,
\bar{v}_0)\nparallel,\]
where $C_4$ is a positive constant
independent of $\|(\bar{\sigma}_0, \bar{v}_0)\nparallel$.
\end{Lemma}

As usual, the global existence will be obtained by a combination of
the local existence result Lemma 5.1 and the a priori estimate below.

\begin{Lemma}
(A priori estimate) Under the assumptions of Lemma 5.1, suppose that
the Cauchy problem (\ref{5.1}), (\ref{5.2}) has a unique classical solution
$(\bar{\sigma},\bar{v})\in \bar{X}(0,T_1)$ for some positive
constant $T_1$. Then there exists two small constants $\delta>0$ and $C_5>0$ which are independent of $T_1$ such that if
\begin{equation}\label{5.5}
\sup_{0\leq t\leq
T_1}\|(\bar{\sigma},\bar{v})(t)\nparallel\leq\delta,
\end{equation}
it holds that
\begin{equation}\label{5.6}
\|\bar{\sigma}(t)\|^2_{N-1}+\|\bar{v}(t)\|^2_{N-2}+\int_0^t\left(\|\nabla\bar{\sigma}(\tau)\|^2_{N-1}+\|\nabla\bar{v}(\tau)\|^2_{N-2}\right)d\tau
\leq
C_5\left(\|\bar{\sigma}_0\|^2_{N-1}+\|\bar{v}_0\|^2_{N-2}\right)
\end{equation}
for all $t\in[0,T_1]$.
\end{Lemma}
\noindent{\bf Proof.}~~Noticing that some smallness
conditions can be imposed on $(\sigma^{per}, v^{per})$, without loss of generality,
we may assume $|||(\sigma^{per}, v^{per})|||\leq\epsilon$ with
$\epsilon>0$ being sufficiently small. Then by the similar argument as in the proof of Lemmas 3.3-3.4, we can obtain
\begin{equation}\label{5.7}
\begin{array}{rl}
 &\displaystyle\frac{d}{dt}\left(\|\bar{U}\|^2+\|\nabla\bar{\sigma}\|^2+d_2\langle
\bar{v}, \nabla\bar{\sigma}\rangle\right)+\|\nabla \bar{v}\|^2+\|\nabla\bar{\sigma}\|^2_1\\[2mm]
 &\leq \epsilon C\left(\|\nabla^3\sigma\|^2_{N-7}+\|\nabla^2\bar{v}\|^2_{N-6}\right),
\end{array}
\end{equation}
and
\begin{equation}\label{5.8}
\begin{array}{rl}
 &\displaystyle\frac{d}{dt}\left(\|\nabla\bar{\sigma}\|^2_{N-2}+\|\nabla \bar{v}\|^2_{N-3}+d_3\sum_{|\alpha|=1}^{N-2}\langle
\partial_x^\alpha \bar{v}, \partial_x^\alpha \nabla\bar{\sigma}\rangle\right)+\|\nabla^2\bar{\sigma}\|^2_{N-2}+\|\nabla^2 \bar{v}\|^2_{N-3}\\[4mm]
 &\leq \epsilon
 C\left(\|\nabla\bar{\sigma}\|^2+\|\nabla\bar{v}\|^2\right),
\end{array}
\end{equation}
 where $d_2>0$ and $d_3>0$ are some suitably small constants, and $C$ is a constant depending only on $\rho_\infty, \mu, \nu$ and $\kappa$. Adding (\ref{5.8})
to (\ref{5.7}), it holds
\begin{equation}\label{5.9}
\begin{array}{rl}
 \displaystyle\frac{d}{dt}& \displaystyle\left(\|\bar{\sigma}\|^2_{N-1}+\|
\bar{v}\|^2_{N-2}+d_2\langle \bar{v},
\nabla\bar{\sigma}\rangle+d_3\sum_{|\alpha|=1}^{N-2}\langle
\partial_x^\alpha \bar{v}, \partial_x^\alpha \nabla\bar{\sigma}\rangle\right)\\[4mm]
 &+\|\nabla\bar{\sigma}\|^2_{N-1}+\|\nabla\bar{v}\|^2_{N-2}\leq 0,
\end{array}
\end{equation}
provided that $\epsilon$ is sufficiently small. Integrating (\ref{5.9}) in $t$ over $(0,t)$, one can immediately get
(\ref{5.6}) since
\[\|\bar{\sigma}\|^2_{N-1}+\|\bar{v}\|^2_{N-2}+d_2\langle \bar{v},
\nabla\bar{\sigma}\rangle+d_3\sum_{|\alpha|=1}^{N-2}\langle
\partial_x^\alpha \bar{v}, \partial_x^\alpha \nabla\bar{\sigma}\rangle\thicksim\|\bar{\sigma}\|^2_{N-1}+\|\nabla
\bar{v}\|^2_{N-2}.\]
by the smallness of $d_2$ and $d_3$. This completes the proof of Lemma 5.2.\\
\noindent{\bf Proof of Theorem 1.2.}~~By Lemmas 5.1-5.2 and the continuity argument,  the Cauchy problem (\ref{5.1}), (\ref{5.2})
admits  a unique solution $(\bar{\sigma}, \bar{v})$ globally in time, which satisfies
(\ref{1.6}) and (\ref{1.7}). Then all the statements in Theorem 1.2 follow immediately. This completes the proof of Theorem 1.2.

\end{document}